\documentclass[a4paper,12pt,leqno]{article}

\usepackage{amsfonts}
\usepackage{latexsym}
\usepackage{epsfig}
\usepackage{amssymb}
\usepackage{graphicx}
\usepackage{oldgerm}
\usepackage{theorem}

\def\R{{\bf R}}
\def\RV{{\bf R}^{N}_{<}}

\def\E0{{\bf E}_{\bf 0}}

\def\x{{\bf x}}
\def\y{{\bf y}}

\def\0{{\bf 0}}

\def\Pf{{\rm Pf}}

\def\cN{\mathcal{N}_N}
\def\cA{\mathcal{A}}
\def\cQ{\mathcal{Q}}
\def\cD{\mathcal{D}}
\def\cS{\mathcal{S}}
\def\cI{\mathcal{I}}
\def\cG{\mathcal{P}}
\def\cO{\mathcal{O}}
\def\a{\mathfrak{a}}
\def\q{\mathfrak{q}}
\def\X{\mathfrak{X}}
\def\bA{\mathbb{A}}
\def\bS{\mathbb{S}}
\def\bD{\mathbb{D}}
\def\bI{\mathbb{I}}
\def\bQ{\mathbb{Q}}
\def\ba{{\bf a}}
\def\bq{{\bf q}}
\def\hq{{\widehat{q}}}
\def\hQ{\widehat{Q}}
\def\tS{\tilde{S}}
\def\tI{\tilde{I}}

\theorembodyfont{\itshape}

\newcommand{\SSC}[1]{\section{#1}\setcounter{equation}{0}}
\begin{document}
\begin{center}
{\bf \Large{Infinite systems of non-colliding Brownian particles}}
\vskip 3mm
{\sc Makoto Katori, Taro Nagao and Hideki Tanemura}
\\
{\it Chuo University, Osaka University} and {\it Chiba University} 
\end{center}

\pagestyle{plain}
\vskip 0.3cm

\begin{small}
\noindent  
{\bf Abstract.} 
Non-colliding Brownian particles in one dimension is studied.
$N$ Brownian particles start from the origin at time $0$
and then they do not collide with each other 
until finite time $T$.
We derive the determinantal expressions for the multitime correlation 
functions using the self-dual quaternion matrices.
We consider the scaling limit of the infinite particles
$N \to \infty$ and the infinite time interval $T \to \infty$.
Depending on the scaling, two limit theorems are
proved for the multitime correlation functions, which may define
temporally inhomogeneous infinite particle systems.
\end{small}

\vspace{1cm}




\normalsize
\SSC{Introduction}

We consider the process $X(t)$, which represents the system of 
$N$ Brownian motions in one dimension all
started from the origin and conditioned never to collide 
with each other up to time $T$.
If we take the limit $T \to \infty$, the system becomes
a temporally homogeneous diffusion process $Y(t)$,
which is the Doob $h$-transform \cite{Do84}
of the absorbing Brownian motion in a Weyl chamber
$$
\RV=\Big\{\x=(x_{1}, x_{2}, \dots, x_{N}) \in \R^{N};
x_{1} < x_{2} < \cdots < x_{N} \Big\},
$$ 
with harmonic function $h_N(\x)= \prod_{1\le i<j\le N}(x_j-x_i)$
\cite{Gra99}. By virtue of the Karlin-McGregor formula 
\cite{KM59_1, KM59_2}, its transition density 
$f_{N}(t, \x, \y)$ from the state $\x$ to $\y$ in $\RV$
in time period $t > 0$ is given by
$$
f_{N}(t, \x, \y)=\det_{1 \leq i, j \leq N} \Big(
p_t (x_{i}, y_{j}) \Big),
$$
where 
$p_t (x,y)= \frac{1}{\sqrt{2 \pi t}} e^{-(x-y)^2/2t}$.
On the other hand, if the non-colliding time interval
$T$ remains finite, the process $X(t), 0 \leq t \leq T$,
is temporally inhomogeneous \cite{KT02b}.

We notice an integral formula found in Harish-Chandra \cite{HC57},
Itzykson and Zuber \cite{IZ80}, and Mehta \cite{Meh81},
$$
\frac{\displaystyle{\det_{1 \leq i, j \leq N}
(p_{t}(x_{i}, y_{j}))}}{h_{N}(\x) h_{N}(\y)}
= c \int dU \ \exp \left[ -\frac{1}{2t}
{\rm tr}(X-U^{\dagger}YU)^2 \right]
$$
with $c^{-1}=(2 \pi)^{N/2} t^{N^2/2} \prod_{i=1}^{N+1} 
\Gamma(i)$, 
where $X$ and $Y$ are the $N \times N$ diagonal matrices,
$X_{ij}=x_{i}\delta_{ij}, Y_{ij}=y_{i}\delta_{ij}$,
and the integral is taken over the group of
unitary matrix $U$ of size $N$. This equality implies that the 
non-colliding Brownian motions such as $X(t)$ and $Y(t)$
can be described by using the eigenvalue-statistics of
Hermitian random matrices in Gaussian ensembles \cite{Meh91}.
In earlier papers\cite{KT02,KT02b}, it was shown that 
$Y(t)$ is identified with Dyson's Brownian motion model
with $\beta=2$ \cite{Dys62} and the particle distribution
is expressed by the probability density of eigenvalues of
random matrices in the Gaussian unitary ensemble (GUE)
with variance $t$, while
$\sqrt{\frac{T}{t(2T-t)}}X(t)$ coincides with the distribution
of eigenvalues of random matrices in the Pandey-Mehta ensemble
\cite{MP83, PM83} with $\alpha = \sqrt{\frac{T-t}{T}}$,
and this temporally inhomogeneous process exhibits
a transition from the GUE statistics
to the Gaussian orthogonal ensemble (GOE) statistics
as the time $t$ goes on from 0 to $T$.

It is known that the eigenvalue distributions of Hermitian
random matrices have determinantal expressions.
For instance, in the GUE, the probability density of
$N$ eigenvalues is expressed by
$$
\rho_{N}(x_{1}, x_{2}, \dots, x_{N}) = \frac{1}{N!}
\det_{1 \leq i, j \leq N} ( K_{N}(x_{i}, x_{j}) ),
$$
with 
$K_{N}(x,y)=\sum_{\ell=0}^{N-1}\varphi_{\ell}(x) \varphi_{\ell}(y)$,
where 
\begin{equation}
\varphi_{\ell}(x)= \frac{1}{\sqrt{h_{\ell}}}
e^{-x^2/2} H_{\ell}(x) 
\label{eqn:varphi}
\end{equation}
with
the $\ell$-th Hermite polynomial $H_{\ell}(x)$ and 
$h_{\ell}=\sqrt{\pi} \, 2^\ell \, \ell !$. By the orthogonality of
$\varphi_{k}(x)$, we can prove the equality
\begin{equation}
\int \det_{1 \leq i, j \leq N'} K_{N}(x_{i}, x_{j}) \, dx_{N'} 
=(N-N'+1) \det_{1 \leq i, j \leq N'-1} K_{N}(x_{i}, x_{j})
\label{eqn:integral}
\end{equation}
for any $1 \leq N' \leq N$. Such integral property enables us
not only to obtain determinantal expressions for correlation
functions, but also to argue the $N \to \infty$ limit
of the system by studying the large $N$ asymptotic of
the function $K_{N}(x,y)$. 
With proper scaling limit, determinantal point processes
with sine-kernel and Airy-kernel are derived.
See \cite{Sos00} and references therein.

In the present paper, we derive the determinantal 
expressions of the multitime correlation functions 
for the process $X(t)$. 
Our aim is to prove limit theorems of the multitime
correlation functions in the scaling limits
of infinite particles $N \to \infty$ and infinite time
interval $T \to \infty$. Depending on the scaling,
we derive two kinds of limit theorems, one of which
provides a spatially homogeneous but temporally
inhomogeneous infinite particle system (Theorem 1), and
other of which does the system with inhomogeneity both in
space and time (Theorem 2). 
We remark that it is easier to prove the limit theorems
for Dyson's Brownian motion model $Y(t)$.
Corresponding to Theorem 1, we will obtain the multitime
correlation functions of the homogeneous system,
which coincides with the system studied by
Spohn \cite{Sp87}, Osada \cite{O96}, and
Nagao and Forrester \cite{NF98b}.
Similarly, corresponding to Theorem 2, an infinite system
with spatial inhomogeneity will be derived, which is
related with the Airy process recently studied by
Pr\"ahofer and Spohn \cite{PS01} and Johansson \cite{Joh02}.

One of the key points of our arguments is that, 
in order to give the determinantal expressions
for the correlation functions for the present processes,
we shall prepare matrices with the elements, which are
neither real nor complex numbers, but quaternions
$$
  q=q_{0}+q_{1} e_{1}+q_{2} e_{2} + q_{3} e_{3} \in {\bf Q}
$$
with $q_{i} \in {\bf C}, 0 \leq i \leq 3$,
in which the four basic units
$\{1, e_{1}, e_{2}, e_{3}\}$ have the following
$2 \times 2$ matrix representations,
$C: {\bf Q} \mapsto {\rm Mat}_{2}({\bf C})$;
\begin{eqnarray}
&& C(1)=\left( 
\begin{array}{cc}
\ 1 \ & \ 0 \ \cr 0 & 1
\end{array}
\right), 
\hskip 2cm
C(e_{1})=\left(
\begin{array}{rr}
\ 0 \ & -1 \cr \ 1 \ & 0
\end{array}
\right), \nonumber\\
&& C(e_{2})=\left( 
\begin{array}{cc}
0 & -\sqrt{-1} \cr -\sqrt{-1} & 0
\end{array}
\right), \quad
C(e_{3})= \left(
\begin{array}{cc}
\sqrt{-1} & 0 \cr 0 & -\sqrt{-1}
\end{array}
\right). \nonumber
\end{eqnarray}
The dual of a quaternion $q$ is defined by
$q^{\dagger}=q_{0}-\sum_{i=1}^{3} q_{i} e_{i}$, and
for a quaternion matrix $Q=(q_{ij}), q_{ij} \in {\bf Q}$,
its dual matrix $Q^{\dagger}=((Q^{\dagger})_{ij})$ is
defined to have the elements 
$(Q^{\dagger})_{ij}=q_{ji}^{\dagger}$.
Following Dyson's definition of the quaternion determinant
for self-dual matrices \cite{D70, Meh89, Meh91}, we can give
the quaternion determinantal expressions having
the similar properties to (\ref{eqn:integral}) for
arbitrary multitime correlation functions 
for $X(t)$ (Theorem 3).
As briefly reported in \cite{NKT}, 
the present results can be regarded as 
simple applications of the results given
in Nagao and Forrester \cite{NF99} and Nagao \cite{N01} 
for multimatrix models, and in
Forrester, Nagao and Honner \cite{FNH99} for
the asymptotic of quaternion determinantal systems,
here we give, however, a self-contained explanation for all the
formulae and calculus developed in the random 
matrix theory, which are used to prove our limit theorems.

The theorems proved here mean the convergence of processes
in the sense of finite dimensional distributions.
As argued in Pr\"ahofer and Spohn \cite{PS01}
and in Johansson \cite{Joh02}, tightness in time should
be confirmed.

\vskip 3mm
\SSC{Statement of Results}

For a given $T > 0$, we define
\begin{equation}
g_{N}^{T}(s, \x; t, \y)
=\frac{f_{N}(t-s, \x, \y) \cN(T-t,\y)}
{\cN(T-s, \x)}
\label{eqn:gNT}
\end{equation}
for $0 \leq s \leq t \leq T, \x, \y \in \RV$,
where
$\cN(t, \x)=\int_{\RV} f_{N}(t, \x, \y) d \y$,
which is the probability that a Brownian motion 
started at $\x\in \RV$ does not hit the boundary of $\RV$
up to time $t>0$.
The function $g_{N}^{T}(s,\x;t,\y)$ can be
regarded as the transition probability density from the
state $\x \in \RV $ at time $s$ 
to the state $\y \in \RV$ at time $t$, 
and associated with the temporally inhomogeneous diffusion process, 
which is the $N$ Brownian motions conditioned not 
to collide with each other in a time interval $[0,T]$.
In \cite{KT02, KT02b} it was shown that
as $|\x|\to 0$, $g_{N}^{T}(0,\x,t,\y)$ converges to 
\begin{equation}
g_{N}^{T}(0, {\bf 0}, t, \y)
\equiv C(N,T,t)  h_{N}(\y) \cN(T-t, \y)
\prod_{i=1}^N p_t(0,y_i),
\label{eqn:gNT0}
\end{equation}
where
$C(N,T,t)=\pi^{N/2}\left(\prod_{j=1}^{N} \Gamma(j/2)\right)^{-1}
T^{N(N-1)/4} t^{-N(N-1)/2}$.
Then the diffusion process $X(t)$ starting from ${\bf 0}$
can be constructed.

We denote by $\X$ the space of countable subset $\xi$ of $\R$
satisfying $\sharp (\xi \cap K) < \infty$
for any compact subset $K$.
We introduce the map $\gamma$ from $\bigcup_{n=1}^{\infty} \R^n$ to
$\X$ defined by
$\gamma (x_1, x_2, \dots, x_n)= \{x_i\}_{i=1}^n$.
Then $\Xi^N (t) = \gamma X(t)$ is the diffusion process on 
the set $\X$ with transition density function
$\widetilde{g}_{N}^{T}(s, \xi ; t, \eta)$, $0 \leq s \leq t \leq T$:
$$
\widetilde{g}_{N}^{T}(s, \xi ; t, \eta)=
\left\{
\begin{array}{ll}
g_{N}^{T}(s, \x; t, \y),
& \mbox{if} \ s>0,  \ \sharp \xi = \sharp \eta =N,
\\
g_{N}^{T}(0, {\bf 0}; t, \y),
& \mbox{if} \ s=0, \ \xi =\{ 0 \}, \ \sharp \eta=N,
\\
0, 
& \mbox{otherwise},
\end{array}\right.
$$
where $\x$ and $\y$ are the elements of $\RV$ with
$\xi=\gamma \x$, $\eta =\gamma\y$.
For $\x^{(m)}_N \in \RV$, $1\leq m \leq M+1$,
and $N'=1,2,\dots, N$, we put 
$\x^{(m)}_{N'} 
= \left(x_1^{(m)}, x_2^{(m)}, \dots, x_{N'}^{(m)}\right)$
and $\xi_m^{N'} = \gamma \x^{(m)}_{N'}$.
For a given time interval $[0,T]$, we consider the $M$
intermediate times
$0 < t_{1} < t_{2} < \cdots < t_{M} < T$.
Then the multitime transition density function 
of the process $\Xi^N (t)$ is given by
\begin{equation}
\label{def:rho}
\rho_N^T(t_{1}, \xi^N_1;\dots; t_{M+1}, \xi^N_{M+1})
=\prod_{m=0}^{M} \widetilde{g}_{N}^{T}(t_{m}, \xi^N_{m}; 
t_{m+1}, \xi^N_{m+1}),
\end{equation}
where, for convenience, we set
$t_{0}=0$, $t_{M+1}=T$ and $\xi^N_0=\{ 0\}$.
From (\ref{eqn:gNT}) and (\ref{eqn:gNT0})
we have 
\begin{eqnarray}
\label{eqn:rho_N}
&& \rho_N^T(t_{1}, \xi^N_1 ; t_2, \xi^N_2 ; \dots ; 
t_{M+1}, \xi^N_{M+1})
= C(N,T,t_{1}) h_{N}\left(\x^{(1)}_N \right) {\rm sgn}
\left( h_{N} \left(\x^{(M+1)}_N \right) \right) \\
&& \hskip 5cm \times
\prod_{i=1}^N p_{t_1}\left(0, x_i^{(1)}\right)
\prod_{m=1}^{M} 
\det_{1 \leq i, j \leq N} \left(
p_{t_{m+1}-t_{m}} \left(x_{i}^{(m)}, x_{j}^{(m+1)} \right) \right).
\nonumber
\end{eqnarray}
For a sequence $\{N_m \}_{m=1}^{M+1}$ of positive integers 
less than or equal to $N$,
we define the \\
$(N_1, N_2,\dots, N_{M+1})$-multitime 
correlation function by
\begin{eqnarray}
\label{def:corr}
&& \rho_N^T \left(t_{1}, \xi^{N_1}_1; t_2, \xi^{N_2}_2 ; 
\dots; t_{M+1}, \xi^{N_{M+1}}_{M+1} \right) 
\\
&& \quad =
\int_{\prod_{m=1}^{M+1} \R^{N-N_{m}}}
\prod_{m=1}^{M+1}
\frac{1}{(N-N_{m})!}\prod_{i=N_{m}+1}^{N} dx_{i}^{(m)}
\rho_N^T(t_{1}, \xi^{N}_1; t_2, \xi^{N}_2, \dots; 
t_{M+1}, \xi^{N}_{M+1}).
\nonumber
\end{eqnarray}

We will study limit theorems of the correlation functions
$\rho^{T_N}_N$ as $N\to \infty$.
First, we consider the case $T_N=2N$.
Let
\begin{eqnarray}
&&\widetilde{\bS}(s, x ; t, y)
\nonumber\\
&&\quad=
\left\{
   \begin{array}{ll}
\displaystyle{
\frac{1}{\pi}
\int_{0}^{1} d \lambda \
\cos(\lambda(x-y)) 
e^{-\lambda^2(t-s)/2}}
 & \mbox{if} \ s > t\\
 & \\
\displaystyle{
\frac{\sin(x-y)}{\pi(x-y )}}
 & \mbox{if} \ s=t \\
 & \\
\displaystyle{
-\frac{1}{\pi} \int_{1}^{\infty} d \lambda \
\cos(\lambda(x-y)) 
e^{-\lambda^2(t-s)/2}}
 & \mbox{if} \ s< t \\
   \end{array}\right. 
\nonumber\\
&&{\bD}(s,x ; t, y)
= -\frac{1}{\pi}
\int_{0}^{1} d \lambda \ \lambda 
\sin(\lambda(x-y))
e^{-(s+t)\lambda^2/2}
\nonumber\\
&&\widetilde{\bI}(s,x ; t, y)
= -\frac{1}{\pi}
\int_{1}^{\infty} d \lambda \
\frac{1}{\lambda} \sin(\lambda(x-y)) 
e^{(s+t)\lambda^2/2}.
\nonumber
\end{eqnarray}
And let ${\bq}^{m,n}(x,y)$ be the quaternion,
whose $2 \times 2$ matrix expression is given by
$$
C({\bq}^{m,n}(x,y))=
\left( 
\begin{array}{cc}
\widetilde{\bS}(s_m,x;s_n,y) & \widetilde{\bI}(s_m,x;s_n,y)\cr
{\bD}(s_m,x;s_n,y) & \widetilde{\bS}(s_n,y;s_m,x)
\end{array}
\right).
$$
Let $M \geq 1$ and $\{N_{m}\}_{m=1}^{M+1}$ be a sequence
of positive integers. 
We denote by 
$\bQ \left(\x^{(1)}_{N_1},\x^{(2)}_{N_2},\dots, 
\x^{(M+1)}_{N_{M+1}} \right)$
the self-dual
$\sum_{m=1}^{M+1} N_m \times \sum_{m=1}^{M+1} N_m$ 
quaternion matrix whose elements are
$\bq^{m,n}\left(x^{(m)}_i, x^{(n)}_j\right)$,
$1\le i\le N_m$, $1\le j\le N_n$,
$1\le m,n \le M+1$, that is,
\begin{eqnarray}
&&\bQ\left(\x^{(1)}_{N_1}, \x^{(2)}_{N_2},\dots, 
\x^{M+1}_{N_{M+1}}\right)
\nonumber
\\
&&=\left[
\begin{array}{ccc}
\bQ^{1,1} \left(\x^{(1)}_{N_1},\x^{(1)}_{N_1}\right) & 
\cdots 
&\bQ^{1, M+1}\left(\x^{(1)}_{N_1},\x^{(M+1)}_{N_{M+1}}\right) 
\cr
\bQ^{2,1}\left(\x^{(2)}_{N_2},\x^{(1)}_{N_1}\right) & 
\cdots 
&\bQ^{2, M+1}\left(\x^{(2)}_{N_2},\x^{(M+1)}_{N_{M+1}}\right) 
\cr
   \cdots    & \cdots &   \cdots     \cr
    \cdots   & \cdots &   \cdots    
\cr
 \bQ^{M+1, 1}\left(\x^{(M+1)}_{N_{M+1}},\x^{(1)}_{N_1}\right) & 
\cdots 
&\bQ^{M+1, M+1}\left(\x^{(M+1)}_{N_{M+1}},\x^{(M+1)}_{N_{M+1}}\right) 
\end{array}
\right]
\nonumber
\end{eqnarray}
with blocks of $N_{m} \times N_{n}$ quaternion matrices
$$
\bQ^{m, n}\left(\x^{(m)}_{N_m}, \x^{(n)}_{N_n}\right)
= \left(
\bq^{m, n}\left(x^{(m)}_i, x^{(n)}_j\right) 
\right)_{1\leq i \leq N_m, 1\leq j \leq N_n},
$$
for $1 \leq m, n \leq M+1$.

For an $N \times N$ self-dual quaternion matrix $Q$, 
the quaternion determinant ${\rm Tdet} Q$ is defined
by Dyson \cite{D70} as
$$
{\rm Tdet} Q=\sum_{\pi \in S_{N}} 
(-1)^{N-\ell(\pi)} \prod_{1}^{\ell(\pi)}
q_{ab}q_{bc} \cdots q_{da},
$$
where $\ell(\pi)$ denotes the number of exclusive
cycles of the form $(a \to b \to c \to \cdots \to d \to a)$
included in a permutation $\pi \in S_{N}$.
\vskip 0.3cm
\noindent{\bf Theorem 1.}{\itshape
\label{thm:main1}
Let $T_N =2N$.
For any $M \geq 1$, 
any sequence $\{N_{m}\}_{m=1}^{M+1}$ of positive integers,
and any strictly increasing sequence $\{ s_m \}_{m=1}^{M+1}$ of
nonpositive numbers with $s_{M+1}=0$,
\begin{eqnarray}
&&\rho\left(s_1, \xi^{N_1}_1 ; s_2, \xi^{N_2}_2 ; \dots ; 
s_M, \xi_M^{N_M}; 0, \xi^{N_{M+1}}_{M+1}\right)
\nonumber\\
&\equiv&
\lim_{N \to \infty} 
\rho^{T_N}_N \left(
T_N+ s_1, \xi^{N_1}_1 ; 
T_N+s_2, \xi^{N_2}_2 ; \dots ; T_N, \xi^{N_{M+1}}_{M+1} \right)
\nonumber\\
&=& {\rm Tdet} \, {\bQ}\left(\x^{(1)}_{N_1},\x^{(2)}_{N_2},
\dots, \x^{(M+1)}_{N_{M+1}} \right).
\nonumber
\end{eqnarray}
}
\noindent{\bf Remark 1.} 
The above system
is spatially homogeneous, since all elements
of the quaternion determinant are functions of
difference of positions, $x_{i}^{(m)}-x_{j}^{(n)}$.
This expresses the bulk property of our infinite particle
system. When $M=1$, the present system is equivalent with
the $N \to \infty$ limit of the two-matrix model
reported by Pandey and Mehta \cite{MP83,PM83}.
In the system defined by Theorem 1, if we
take the further limit such that
$s_m \to -\infty$ with the time difference
$s_n -s_m$ fixed, $1 \leq m, n \leq M$,
then $\bD (s_m,x;s_n,y) \to \infty$, 
$\widetilde{\bI}(s_m,x;s_n,y)\to 0$, while the product
$\bD(s_m,x;s_n,y) \widetilde{\bI}(s_m,x;s_n,y)\to 0$.
Therefore, we may replace $\bD$ and $\widetilde{\bI}$ by
zeros in this limit, and the quaternion determinant 
${\rm Tdet} \ {\bQ}\left(\x^{(1)}_{N_1},\x^{(2)}_{N_2},
\dots, \x^{(M)}_{N_{M}} \right)$
will be reduced to an ordinary determinant
$\det \, {\bA}\left(\x^{(1)}_{N_1},\x^{(2)}_{N_2},\dots, 
\x^{(M)}_{N_{M}} \right)$ with the elements
$\ba^{m,n}\left(x^{(m)}_i,x^{(n)}_j\right)
=\widetilde{\bS}\left(s_m,x^{(m)}_i;s_n,x^{(n)}_j\right)$.
Hence, we obtain a temporally and spatially
homogeneous system, whose correlation functions are given by
$$
{\rho}'
\left(s_1, \xi^{N_1}_1 ; s_2, \xi^{N_2}_2 ; \dots ; 
s_M, \xi_M^{N_M} \right)
= \det {\bA} \left(\x^{(1)}_{N_1},\x^{(2)}_{N_2},
\dots, \x^{(M)}_{N_{M}} \right).
$$
Such a homogeneous system was studied by 
Spohn \cite{Sp87}, Osada \cite{O96} and 
Nagao and Forrester \cite{NF98b}
as an infinite particle limit of 
Dyson's Brownian motion model\cite{Dys62}.

Next, we consider the case that $T_N = 2N^{1/3}$.
In order to state the result, we have to introduce
the following functions.
Let ${\rm Ai}(z)$ be the Airy function:
\begin{equation}
  {\rm Ai}(z)= \frac{1}{2\pi} \int_{-\infty}^{\infty}
e^{\sqrt{-1}(z t + t^3/3)} dt.
\label{eqn:Ai}
\end{equation}
For $s, t<0$ and $x,y \in \R$,
we put
\begin{eqnarray}
{\cD}(s, x ; t, y) 
&=&\frac{1}{4} \left[ \int_{0}^{\infty} d\lambda \ 
e^{s \lambda/2} {\rm Ai}(x+\lambda)
\frac{d}{d\lambda} \left\{ e^{t \lambda/2} {\rm Ai}(y+\lambda) 
\right\} \right. \nonumber\\
&& \left. - \int_{0}^{\infty} d\lambda \ 
e^{t \lambda/2} {\rm Ai}(y+\lambda)
\frac{d}{d\lambda} 
\left\{ e^{s \lambda/2} {\rm Ai}(x +\lambda) \right\} \right], 
\nonumber\\
\widetilde{{\cI}}(s, x ; t, y)
&=& \int_{0}^{\infty} d\lambda \ e^{t \lambda/2} {\rm Ai}(y-\lambda)
\int_{\lambda}^{\infty} d\lambda' \ e^{s \lambda'/2} {\rm Ai}(x-\lambda') 
\nonumber\\
&& - \int_{0}^{\infty} d\lambda \ e^{s \lambda/2} {\rm Ai}(x-\lambda)
\int_{\lambda}^{\infty} d\lambda' \ e^{t \lambda'/2} {\rm Ai}(y-\lambda'),
\nonumber
\end{eqnarray}
and
$$
\widetilde{\cS}(s,x;t,y)=\cS(s,x;t,y)-{\cG}(s,x;t,y) 1(s<t),
$$
with
\begin{eqnarray}
{\cS}(s, x ; t, y) 
&=& \int_{0}^{\infty} d\lambda \ e^{(t -s)\lambda/2}
{\rm Ai}(x+\lambda) {\rm Ai}(y+\lambda) 
+ \frac{1}{2} {\rm Ai}(y)
\int_{0}^{\infty} d\lambda \ e^{s \lambda/2} {\rm Ai}(x-\lambda),
\nonumber\\
{\cG}(s,x ; t,y)
&=& \int_{-\infty}^{\infty} d\lambda \
e^{(t - s)\lambda/2} {\rm Ai}(x+\lambda){\rm Ai}(y+\lambda),
\nonumber
\end{eqnarray}
where $1 (s<t)=1$ if $s<t$, and $=0$ otherwise.
And let ${\q}^{m,n}(x,y)$ be the quaternion,
whose $2 \times 2$ matrix expression is given by
$$
C(\q^{m,n}(x,y))=
\left( 
\begin{array}{cc}
\widetilde{\cS}(s_m,x;s_n,y) & \widetilde{\cI}(s_m,x;s_n,y) \cr
{\cD}(s_m,x;s_n,y) & \widetilde{\cS}(s_n,y;s_m,x) 
\end{array}
\right).
$$
Let $M \geq 1$ and $\{N_{m}\}_{m=1}^{M+1}$ be a sequence
of positive integers. 
We denote by
$\cQ \left(\x^{(1)}_{N_1},\x^{(2)}_{N_2},\dots, 
\x^{(M+1)}_{N_{M+1}} \right)$
the self-dual
$\sum_{m=1}^{M+1} N_m \times \sum_{m=1}^{M+1} N_m$ 
quaternion matrix whose elements are
$\q^{m,n}\left(x^{(m)}_i, x^{(n)}_j\right)$,
$1\le i\le N_m$, $1\le j\le N_n$,
$1\le m,n \le M+1$.
\vskip 3mm
\noindent{\bf Theorem 2.}\
{\itshape
\label{thm:main2}
Let $T_N = 2N^{1/3}$ and $a_N(s) = 2N^{2/3}- s^2/4$ for $s\in \R$.
For any $M \geq 1$, 
any sequence $\{N_{m}\}_{m=1}^{M+1}$ of positive integers,
and any strictly increasing sequence $\{ s_m \}_{m=1}^{M+1}$ of
nonpositive numbers with $s_{M+1}=0$,
\begin{eqnarray}
&&\widehat{\rho} \left(s_1, \xi^{N_1}_1 ; 
\dots ; s_{M+1}, \xi^{N_{M+1}}_{M+1}\right)
\nonumber\\
&\equiv&
\lim_{N \to \infty} 
\rho^{T_N}_N \left(
T_N+ s_1, \theta_{a_N(s_1)}\xi^{N_1}_1 ; 
\dots ;T_N, \theta_{ a_N(s_{M+1})} \xi^{N_{M+1}}_{M+1} \right)
\nonumber\\
&=& {\rm Tdet} \, {\cQ}\left(\x^{(1)}_{N_1},
\dots, \x^{(M+1)}_{N_{M+1}}\right),
\nonumber
\end{eqnarray}
where $\theta_u \{x_i\} = \{ x_i+u \}$.
}
\vskip 3mm
\noindent
{\bf Remark 2.}
This theorem may define an infinite particle system, in which
any type of space-time correlation function is given by
the above quaternion determinant.
This quaternion determinantal system is the same
as that derived in Forrester, Nagao and Honner \cite{FNH99},
and it is inhomogeneous both in space and time.
The spatial inhomogeneity is attributed to the fact that
this system expresses the edge property of the infinite
non-colliding Brownian particles. Thus, if we take the
bulk limit, $x_{i}^{(m)} \to - \infty$ with the position
differences $x_{i}^{(m)}-x_{j}^{(n)}$ fixed,
then the system should recover spatial homogeneity.
It is confirmed by observing that the quaternion determinantal
system given in Theorem 1 can be derived as the bulk limit
of the system of Theorem 2, if we use the asymptotic expansion
of the Airy function (\ref{eqn:Ai}) \cite{AS65},
$$
{\rm Ai}(-x) \sim \frac{1}{\pi^{1/2} x^{1/4}}
\cos \left( \frac{2}{3} x^{3/2}-\frac{\pi}{4} \right)
\quad \mbox{as} \ x \to \infty.
$$
On the other hand, keeping the spatial inhomogeneity, one can
consider the limit
$s_m \to -\infty$ with the time difference
$s_n - s_m$ fixed, $1 \leq m, n \leq M$.
In this limit, 
${\cD}(s_m,x;s_n,y) \to 0$, $\widetilde{\cI}(s_m,x;s_n,y) \to 0$,
and
$$
{\cS}(s_m,x;s_n,y) \to
\int_{0}^{\infty} d\lambda \ e^{(s_n - s_m)\lambda/2}
{\rm Ai}(x+\lambda) {\rm Ai}(y+\lambda).
$$
Hence the off-diagonal elements vanish in the $2 \times 2$
matrix expressions of quaternion ${\q}^{m,n}(x,y)$ and
$$
C(\q^{m,n}(x,y))\to
\left( 
\begin{array}{cc}
\a^{m,n}(x,y) & 0 \cr
0 & \a^{n,m}(y,x) 
\end{array}
\right)
$$
for $1 \leq m,n \leq M$, where
\begin{eqnarray}
&&\a^{m,n}(x,y)=
\a(s_m, x; s_n, y) 
\nonumber
\\
&&=
\left\{
   \begin{array}{ll}
      \displaystyle{
\int_{0}^{\infty} d\lambda \ e^{(s_n -s_m)\lambda/2}
{\rm Ai}(x+\lambda) {\rm Ai}(y+\lambda) }
 & \mbox{if} \ m \geq n \\
      \displaystyle{
-\int_{-\infty}^{0} d\lambda \ e^{(s_n -s_m)\lambda/2}
{\rm Ai}(x+\lambda) {\rm Ai}(y+\lambda)}
 & \mbox{if} \ m < n. \\
   \end{array}\right.
\nonumber
\end{eqnarray}
Then the quaternion determinant
${\rm Tdet} {\cQ}\left(\x^{(1)}_{N_1},\x^{(2)}_{N_2},
\dots, \x^{(M)}_{N_{M}} \right)$ 
is reduced to an ordinary determinant 
$\det \cA \left(\x^{(1)}_{N_1},\x^{(2)}_{N_2},
\dots, \x^{(M)}_{N_{M}} \right)$ with the elements
$\a^{m,n}\left(x^{(m)}_i, x^{(n)}_j \right)$.
In this way, we will obtain the infinite particle system, which
is temporally homogeneous but spatially inhomogeneous
with the multitime correlation functions
$$
\widehat{\rho}'
\left(s_1, \xi^{N_1}_1 ; s_2, \xi^{N_2}_2 ; \dots ; 
s_M, \xi_M^{N_M}\right)
= \det \, {\cA}\left(\x^{(1)}_{N_1},\x^{(2)}_{N_2},
\dots, \x^{(M)}_{N_{M}} \right).
$$
In particular, if we set
$N_{1}=N_{2}=\cdots=N_{M}=1$,
then
$$
\widehat{\rho}'\left(s_1, \{x^{(1)} \}; \dots,
s_{M}, \{x^{(M)} \}\right)
=\det_{1 \leq m,n \leq M} \left(
\a^{m,n}\left(x^{(m)}, x^{(n)}\right)
\right).
$$
This is the same as the system called
the Airy process by Pr\"ahofer and Spohn in \cite{PS01}.
(See also \cite{Joh02}.)

\SSC{Quaternion determinantal expressions of the correlations}

In this section we give quaternion determinantal expressions for 
the correlation functions defined in (\ref{def:corr})
along the procedure in \cite{N01}.
From now on we consider the case $N$ is even, 
for simplicity of notations. See \cite{N01}, for
necessary modifications for odd case.
For $1 \leq m, n \leq M+1$, define
\begin{equation}
F^{m,n}(x,y)
= \int_{-\infty}^{\infty} dw \int_{-\infty}^{w} dz \ 
\left| \begin{array}{cc}
 p_{T-t_m}(x,z) & p_{T-t_m}(x, w)
 \cr
 p_{T-t_n}(y, z) & p_{T- t_n}(y, w)
\end{array}
\right|,
\label{def:Fmn}
\end{equation}
where $p_0 (y,x)dy = \delta_x (dy)$.
We introduce an antisymmetric inner products
$$
\langle f, g \rangle_{m} =
\int_{\R} dx \int_{\R} dy \ F^{m,m}(x,y) f(x) g(y),
$$
and
$$
\langle f, g \rangle =
\int_{\R}  dx \int_{\R} dy \ F^{1,1}(x,y)
p_{t_1}(0,x)p_{t_1}(0,y) f(x) g(y).
$$
For $k=0,1,\dots$ we consider the polynomials in $x$
of degree $k$ defined by
\begin{equation}
\label{def:Rk}
R_{k}(x) = z_1^{-k}
\sum_{j=1}^k \alpha_{k j}H_j\left( \frac{x}{c_1}\right)
z_1^{j},
\end{equation}
where
$c_1 = \sqrt{ \frac{t_1(2T-t_1)}{T}}$,
$z_1 = \sqrt{\frac{2T-t_1}{t_1}}$,
\begin{equation}
\alpha_{k j}=\left\{
\begin{array}{ll}
2^{-k}c_1^{k}\delta_{k j},
& \mbox{if $k$ is even,}
\\ 
2^{-k}c_1^{k}
\Big\{ \delta_{k j}- 2(k-1)\delta_{k-2 \, j} \Big\},
& \mbox{if $k$ is odd,}
\end{array}
\right.
\label{def:alpha}
\end{equation}
and $H_j(x)$ are the Hermite polynomials.
They are monic and satisfy the skew orthogonal relations:
\begin{eqnarray}
&& \langle R_{2j}, R_{2\ell+1} \rangle
= - \langle R_{2\ell+1}, R_{2j} \rangle
= r_{j} \delta_{j \ell}, \nonumber\\
&& \langle R_{2j}, R_{2 \ell} \rangle
= \langle R_{2j+1}, R_{2 \ell+1} \rangle
=0, \quad j, \ell=0,1,2, \cdots,
\nonumber
\label{eqn:skew01}
\end{eqnarray}
where
$$
r_j = \frac{\Gamma (j+\frac{1}{2})\Gamma (j+1)}{\pi} 
\left( \frac{t_{1}^2}{T} \right)^{2j+1/2}.
$$
For $m=1,2,\dots, M+1$, and $k=0,1,\dots$, put
\begin{equation}
R^{(m)}_{k}(x) = \int_{\R} dy \ R_{k}(y) p_{t_1}(0,y)p_{t_m-t_1}(y,x).
\label{def:Rkm}
\end{equation}
Then we can prove the skew orthogonal relations
\begin{eqnarray}
&& \langle R_{2j}^{(m)}, R_{2\ell+1}^{(m)}\rangle_{m}
= - \langle R_{2\ell+1}^{(m)}, R_{2j}^{(m)} \rangle_{m}
= r_{j} \delta_{j \ell}, 
\nonumber\\
&& \langle R_{2j}^{(m)}, R_{2 \ell}^{(m)} \rangle_{m}
=\langle R_{2j+1}^{(m)}, R_{2 \ell+1}^{(m)} \rangle_{m} =0,
\quad j, \ell=0,1,2, \cdots,
\nonumber
\label{eqn:skew02}
\end{eqnarray}
for any $m=1,2,\dots, M+1$.
For $m=1,2,\dots, M+1$, define
\begin{equation}
\Phi^{(m)}_{k}(x) = \int_{\R} dy \ R^{(m)}_{k}(y)
F^{m,m}(y, x), \quad
k=0,1,2, \dots. 
\label{def:Phi}
\end{equation}
Now we introduce the functions on $\R^2$,
$D^{m,n}, I^{m,n}$ and $S^{m,n}, 1 \leq m,n \leq M+1$,
given by
\begin{eqnarray}
&&D^{m,n}(x,y)
= \sum_{k=0}^{(N/2)-1} 
\frac{1}{r_{k}}\Big[R_{2k}^{(m)}(x) R_{2k+1}^{(n)}(y) 
- R_{2k+1}^{(m)}(x) R_{2k}^{(n)}(y)\Big],
\label{def:Dmn}
\\
&&I^{m,n}(x,y)
= - \sum_{k=0}^{(N/2)-1} 
\frac{1}{r_{k}}\Big[\Phi_{2k}^{(m)}(x) \Phi_{2k+1}^{(n)}(y)
- \Phi_{2k+1}^{(m)}(x) \Phi_{2k}^{(n)}(y)\Big],
\label{def:Imn}
\\
&&S^{m,n}(x,y)
= \sum_{k=0}^{(N/2)-1} 
\frac{1}{r_{k}}\Big[\Phi_{2k}^{(m)}(x) R_{2k+1}^{(n)}(y) 
- \Phi_{2k+1}^{(m)}(x) R_{2k}^{(n)}(y)\Big].
\label{def:Smn}
\end{eqnarray}
Further we define
\begin{eqnarray}
&&\tS^{m,n}(x,y)=S^{m,n}(x,y)- p_{t_n-t_m}(x,y)1(m < n),
\label{def:tSmn}
\\
&&\tI^{m,n}(x,y) = I^{m,n}(x,y) + F^{m,n}(x, y).
\label{def:tImn}
\end{eqnarray}
Define the quaternions 
$q^{m, n}(x,y), 1 \leq m, n \leq M+1, x,y \in \R$ 
so that these $2 \times 2$ matrix expressions 
$C(q^{m,n}(x,y))$ are given by
$$
C(q^{m, n}(x,y))=\left( 
\begin{array}{cc}
\widetilde{S}^{m, n}(x,y) & \widetilde{I}^{m, n}(x,y) \\ 
D^{m, n}(x,y) & \widetilde{S}^{n, m}(y,x) 
\end{array}
\right).
$$

Let $M \geq 1$ and $\{N_{m}\}_{m=1}^{M+1}$ be a sequence
of positive integers less than or equal to $N$. 
For $\x^{(m)}_N \in \RV$, $1 \leq m \leq M+1$,
we denote by 
$Q \left(\x^{(1)}_{N_1},\x^{(2)}_{N_2},\dots, 
\x^{(M+1)}_{N_{M+1}} \right)$
the self-dual
$\sum_{m=1}^{M+1} N_m \times \sum_{m=1}^{M+1} N_m$ 
quaternion matrix whose elements are
$q^{m,n}\left(x^{(m)}_i, x^{(n)}_j \right)$,
$1\le i\le N_m$, $1\le j\le N_n$,
$1\le m,n \le M+1$.
Then we show the following relation.

\vskip 3mm
\noindent{\bf Theorem 3.}\
{\itshape
The multitime correlation function {\rm (\ref{def:corr})}
is written as
$$
\rho^T_N \left(t_{1}, \xi^{N_1}_1; \dots; 
t_{M+1}, \xi^{N_{M+1}}_{M+1} \right)
= {\rm Tdet} Q \left(\x^{(1)}_{N_1},\x^{(2)}_{N_2},\dots, 
\x^{(M+1)}_{N_{M+1}} \right).
$$
}

In order to prove the theorem,
first we introduce the Pfaffian.
For an integer $N$ and an antisymmetric $2N \times 2N$ matrix
$A=(a_{ij})$, the Pfaffian is defined as
\begin{eqnarray}
&&\Pf(A) = \Pf_{1 \leq i < j \leq 2N}(a_{ij}) \nonumber\\
&&= \frac{1}{N !} \sum_{\sigma}
{\rm sgn}(\sigma) a_{\sigma(1) \sigma(2)} a_{\sigma(3) \sigma(4)} 
\cdots a_{\sigma(2N-1) \sigma(2N)},
\nonumber
\end{eqnarray}
where the summation is extended over all permutations $\sigma$
of $(1,2,\dots, 2N)$ with restriction
$\sigma(2k-1) < \sigma(2k), k=1,2,\dots, N$.
If $Q$ is an $N \times N$ self-dual quaternion matrix, then
\begin{equation}
\label{eqn:Td=Pf}
{\rm Tdet} Q = \Pf \Big( J C(Q) \Big),
\end{equation}
where $J$ is an $2N \times 2N$ antisymmetric matrix with only
non-zero elements
$$
J_{2k+1, 2k+2}=-J_{2k+2,2k+1}=1,
\quad k=0,1,2, \dots, N-1.
$$
See, for instance, Mehta \cite{Meh89}.

For a function $\psi^{m,n}$ defined on $\R^2$
we denote the $N\times N$-matrices whose
$(i,j)$-entry is 
$\psi^{m,n} \left(x_{i}^{(m)}, x_{j}^{(n)} \right)$
by $\psi^{m,n} \left(\x^{(m)}_N, \x^{(n)}_N \right)$,
or simply by $\psi^{m,n}$ for short.
And we denote by
$R^{(m)}\left(\x^{(m)}_N \right)$ the $N\times N$ matrix with
$R^{(m)}\left(\x^{(m)}_N \right)_{i,j}= R^{(m)}_{j-1}(x_i)$,
and by $\Phi^{(m)}\left(\x^{(m)}_N \right)$  that with
$\Phi^{(m)}\left(\x^{(m)}_N \right)_{i,j}= \Phi^{(m)}_{j-1}(x_i)$.
Let $L$ be the $N\times N$ diagonal matrix
with 
$L_{i,i}= \sqrt{\, r_{[(i-1)/2]}}$, $i=1,2,\dots,N$,
and 
$\widetilde{R}^{(m)}\left(\x^{(m)}_N \right)
= L^{-1}R^{(m)}\left(\x^{(m)}_N \right)$.
Then we have 
\begin{equation}
\label{eqn:RJR}
\widetilde{R}^{(m)}\left(\x^{(1)}_N \right) J 
\widetilde{R}^{(n)}\left(\x^{(1)}_N \right)^T
= D^{m,n}\left(\x^{(m)}_N,\x^{(n)}_N \right).
\end{equation}

As the first step of the proof of Theorem 3.
We show that the multitime probability density 
defined in (\ref{def:rho}) is written as
\begin{equation}
\label{eq:pd=TdQ}
\rho^T_N\left(t_1,\xi^N_1;\dots ;
t_{M+1},\xi^N_{M+1} \right)
={\rm Tdet} Q\left(\x^{(1)}_N,\dots,\x^{(M+1)}_N \right).
\end{equation}
For simplicity of notation,
here we give the proof of (\ref{eq:pd=TdQ}) for $M=2$.
It is straightforward to prove (\ref{eq:pd=TdQ}) for general $M$. 
Since
$$
{\rm sgn}\left( h_{N} \left(\x^{(3)}_N \right) \right)
= \Pf_{1 \leq i < j \leq N} \left(
{\rm sgn} \left(x_{j}^{(3)}-x_{i}^{(3)} \right) \right),
$$
and
${\rm sgn}(y-x)=F^{3,3}(x,y)$,
we have
\begin{equation}
\label{eqn:sgnh}
{\rm sgn}\left( h_{N}\left(\x^{(3)}_N \right) \right)
= \Pf \left[F^{3,3} \right].
\end{equation}
Noting that $R_{k}(x)$ is the monic polynomial of degree $k$,
we have
$$
h_{N}\left(\x^{(1)} \right) = \det_{1 \leq i, j \leq N}
\left( \left(x_{j}^{(1)} \right)^{i-1} \right) 
= \det_{1 \leq i, j \leq N} 
\left( R_{i-1}\left(x_{j}^{(1)}\right) \right),
$$
and so
\begin{equation}
\label{eqn:phn=dR1}
\prod_{i=1}^N p_{t_1}\left(0,x^{(1)}_i \right)
h_{N} \left(\x^{(1)} \right) 
= \det\left [ R^{(1)} \left(\x^{(1)}_N \right) \right]
\end{equation}
Since $\det L =\prod_{k=0}^{N/2-1} r_{k} = C(N,T,t_{1})^{-1}$,
from (\ref{eqn:RJR}) and (\ref{eqn:phn=dR1})
\begin{eqnarray}
\label{eqn:Cph=D}
&&C(N,T,t_{1})\prod_{i=1}^N 
p_{t_1} \left(0,x^{(1)}_i \right) h_{N} \left(\x^{(1)} \right) 
= \det \left[ \widetilde{R}^{(1)} \left(\x^{(1)}_N \right) \right] 
\\
&&\qquad\qquad \qquad \qquad = \Pf \
\left[ \widetilde{R}^{(1)} \left(\x^{(1)}_N \right)
 J \widetilde{R}^{(1)}\left(\x^{(1)}_N \right)^T \right]
\nonumber\\
&& \qquad\qquad \qquad \qquad =\Pf 
\left[ D^{1,1}\left(\x^{(1)}_N, \x^{(1)}_N \right) \right].
\nonumber
\end{eqnarray}
Then from (\ref{eqn:rho_N}), (\ref{eqn:sgnh}) and (\ref{eqn:Cph=D})
we have
\begin{eqnarray}
&& \rho_N^T(t_{1}, \xi^N_1; t_2, \xi^N_2 ; 
t_{3}, \xi^N_{3})
\nonumber\\
&&= \Pf[D^{1,1}] \Pf[F^{3,3}]
\prod_{m=1}^2 \det_{1 \leq i, j \leq N} 
\left[p_{t_{m+1}-t_{m}} \right]
\nonumber\\
&&=(-1)^{3N/2}\Pf\left[ 
\begin{array}{cc}
D^{1,1} & O \cr O & -F^{3,3}
\end{array} \right]
\prod_{m=1}^2 
\Pf\left[ 
\begin{array}{cc}
O & -(p_{t_{m+1}-t_{m}})^T \cr p_{t_{m+1}-t_{m}} & O
\end{array} \right].
\nonumber
\end{eqnarray}
By basic properties of the Pfaffians, we have
\begin{eqnarray}
&&\Pf\left[ 
\begin{array}{cc}
D^{1,1} & O \cr O & -F^{3,3}
\end{array} \right]
\prod_{m=1}^2 
\Pf\left[ 
\begin{array}{cc}
O & -(p_{t_{m+1}-t_{m}})^T \cr p_{t_{m+1}-t_{m}} & O
\end{array} \right]
\nonumber\\
&&=\Pf \left[
\begin{array}{cccccc}
D^{1,1} & O             &  O      &  O              &  O      & O \cr
 O      & -F^{3,3} &  O      &  O              &  O      & O \cr
 O      &  O            &  O      & - (p_{t_{2}-t_{1}})^{T} & O & O \cr
 O      &  O            & p_{t_{2}-t_{1}} &  O &  O      & O \cr
 O      &  O            &  O      &  O    & O & -(p_{t_{3}-t_{2}})^{T} \cr
 O      &  O            &  O      &  O    & p_{t_{3}-t_{2}}& O    
\nonumber       
\end{array}
\right] 
\nonumber\\
&&=\Pf \left[
\begin{array}{cccccc}
D^{1,1} & O            &  O       &  O              & O      & O \cr
 O      & O            &  p_{t_{2}-t_{1}} &  O      & O      & O \cr
 O      & -(p_{t_{2}-t_{1}})^T &  O       &  O      & O      & O \cr
 O      &  O           &  O       &  O    & p_{t_{3}-t_{2}}& O \cr
 O      &  O           &  O       &  -(p_{t_{3}-t_{2}})^{T} & O & O \cr
 O      &  O           &  O       &  O              & O      &-F^{3,3}    
\nonumber       
\end{array}
\right] 
\nonumber\\
&&=\Pf \left[
\begin{array}{cccccc}
D^{1,1} & O           & O       & O             & O & O \cr
    &     &     &   &   &  \cr
 O  &-F^{1,1} & p_{t_{2}-t_{1}} &-F^{1,2} & p_{t_{3}-t_{1}} &-F^{1,3} \cr
     &     &     &   &   &  \cr
 O      &-(p_{t_{2}-t_{1}})^T & O       & O     & O & O \cr
    &     &     &   &   &  \cr
 O  &-F^{2,1}     & O       &-F^{2,2}   & p_{t_{3}-t_{2}} &-F^{2,3} \cr
     &     &     &   &   &  \cr
 O  &-(p_{t_{3}-t_{1}})^T & O   &-(p_{t_{3}-t_{2}})^{T} & O & O \cr
    &     &     &   &   &  \cr
 O      &-F^{3,1}     & O       &-F^{3,2}       & O       &-F^{3,3}    
\nonumber       
\end{array}
\right].
\nonumber
\end{eqnarray}
Since $\x^{(1)}_N \in \RV$, $h_N \left(\x^{(1)}_N \right)\not= 0$, 
and so $\det \left[R^{(1)} \left(\x^{(1)}_N \right) \right] 
\not= 0$ by (\ref{eqn:phn=dR1}).
Hence we can define matrices
$$
U^{(m)}=R^{(m)}\left(\x^{(m)}_N \right) 
R^{(1)}\left(\x^{(1)}_N \right)^{-1}, \quad
V^{(m)}= \Phi^{(m)} \left(\x^{(m)}_N \right) 
R^{(1)}\left(\x^{(1)}_N \right)^{-1}, 
$$
which satisfies
\begin{eqnarray}
&U^{(m)} D^{1,1} (U^{(n)})^T = D^{m,n},
\qquad
V^{(m)} D^{1,1} (V^{(n)})^T = -I^{m,n},
\nonumber
\\
&V^{(m)} D^{1,1} (U^{(n)})^T = S^{m,n},
\qquad
U^{(m)} D^{1,1} (V^{(n)})^T = - (S^{n,m})^T.
\nonumber
\end{eqnarray}
By repeating elementary operations, we see that 
the last Pfaffian equals to
\begin{eqnarray}
&&\Pf \left[
\begin{array}{cccccc}
 D^{1,1} &(S^{1,1})^T   & D^{1,2}  &(S^{2,1})^T & D^{1,3} & (S^{3,1})^T \cr
-S^{1,1} &-\tI^{1,1}    &-\tS^{1,2} &-\tI^{1,2} &-\tS^{1,3} &-\tI^{1,3} \cr
 D^{2,1} &(\tS^{1,2})^T & D^{2,2}  & (S^{2,2})^T & D^{2,3}   & S^{3,2} \cr
-S^{2,1} &-\tI^{2,1}    &-\tS^{2,2} &-\tI^{2,2}  &-\tS^{2,3} &-\tI^{2,3} \cr
 D^{3,1} &(\tS^{1,3})^T & D^{3,2}  &(\tS^{2,3})^T & D^{3,3} & (S^{3,3})^T \cr
-S^{3,1} &-\tI^{3,1}    &-S^{3,2}   &-\tI^{3,2} &-S^{3,3}   &-\tI^{3,3}    
\nonumber
\end{array}
\right]
\nonumber\\
&&
=(-1)^{3N/2}
\Pf \left[
\begin{array}{ccc}
A^{1,1} & A^{1,2} & A^{1, 3} \cr
A^{2,1} & A^{2,2} & A^{2,3} \cr
A^{3,1} & A^{3,2} & A^{3,3} 
\nonumber
\end{array}
\right],
\nonumber
\end{eqnarray}
where each $A^{m,n}=(A^{m,n}_{ij})$ is a
$2N \times 2N$ matrix which consists of $2 \times 2$
blocks
$$
A^{m,n}_{ij} = 
\left( 
\begin{array}{cc}
D^{m,n}_{ij} & \tS^{n,m}_{ji} \cr
   &    \cr
-\tS^{m,n}_{ij} & - \tI^{m,n}_{ij} 
\nonumber
\end{array}
\right).
$$
We can see that the above matrix $A=(A^{m,n}_{ij})$ 
satisfies the relation $A=J C(Q)$.
Therefore, (\ref{eq:pd=TdQ}) is derived from (\ref{eqn:Td=Pf}).

For square integrable functions $\phi$ and $\psi$ defined on $\R^2$,
put $\phi * \psi(x,y) = \int_{\R} \phi (x,z) \psi (z, y) dz$.
Then we have
\begin{eqnarray}
&&S^{m,p}*S^{p,m} =I^{m,p}* D^{p,n} =D^{m,p}*F^{p,n}= S^{m,p},
\nonumber\\
&&D^{m,p}* S^{p,n}= D^{m,n},
\quad S^{m,p}* I^{p,n} = S^{m,p}*F^{p,n} = I^{m,n},
\nonumber\\
&&S^{m,p}* p_{t_n-t_p} = S^{m,n}, \quad 
D^{m,p}* p_{t_n-t_p} = D^{m,n}, \quad \mbox{ if $p <n$}.
\nonumber
\end{eqnarray}
Hence by simple calculation we see that
\begin{eqnarray}
\int_{\R} q^{m,m}(z,z)dz &=& N,
\nonumber\\
\int_{\R} q^{m,p}(x,z) q^{p,n}(z,y) dz &=& q^{m,n}(x,y)
+q^{m,n}(x,y)\kappa (n,p) - \kappa (p,m)q^{m,n}(x,y),
\nonumber
\end{eqnarray}
where $\kappa (n,p)$ is a quaternion
with
$$
C(\kappa(n,p)) = \left( 
\begin{array}{cc}
1- 1(p < n) & 0 \cr
\cr 0 & - 1(n <p)
\end{array}
\right).
$$
Then by slight modification of Theorem 6 in \cite{NF99}
we have the following integral formula for any 
$ 1 \leq N_{m} \leq N, m=1,2,\dots, M+1$,
\begin{eqnarray}
&&\int_{\R}{\rm Tdet}
Q \left(\x^{(1)}_{N_{1}},\dots, \x^{(m)}_{N_{m}},\dots, 
\x^{(M+1)}_{N_{M+1}} \right)
dx^{(m)}_{N_{m}}
\nonumber\\
&&\qquad \qquad =(N-N_{m}+1)
{\rm Tdet}Q \left(\x^{(1)}_{N_{1}},\dots, 
\x^{(m)}_{N_{m}-1},\dots, 
\x^{(M+1)}_{N_{M+1}} \right),
\nonumber
\end{eqnarray}
which is the generalization of the formula (\ref{eqn:integral})
given in Introduction of the present paper.
Successive application of the above relation yields Theorem 3.

\vskip 3mm

\SSC{Expansion using Hermite polynomials}

In this section we show expansions of functions $p_{t_n-t_m}$,
$R^{(m)}_k$ and $\Phi^{(m)}_k$ by using Hermite polynomials
$H_k$.
Put 
$$
c_n =\displaystyle{ \sqrt{ \frac{t_n(2T-t_n)}{T} } },
\quad 
\gamma_{n} = \displaystyle{-\frac{T-t_{n}}{T}},
\quad
z_n = \displaystyle{ \sqrt{ \frac{2T-t_n}{t_n} } },
$$
and $\tau^{(n)} = -\log z_n$.
By simple calculation we have
\begin{eqnarray}
p_{t_n -t_m}(x,y)
&=& \frac{e^{- (t_{m}/2T) (x/c_{m})^2 }
e^{ (t_{n}/2T)( y/c_{n})^2 }}{\sqrt{2\pi (t_{n}-t_{m})}}
\nonumber\\
&&\times \exp \left( - \frac{\left\{(y/c_{n})
-e^{-(\tau^{(n)}-\tau^{(m)})} (x/c_{m})\right\}^2}
{1-e^{-2( \tau^{(n)}-\tau^{(m)})})} \right)
\nonumber
\end{eqnarray}
for $1 \leq m < n < M+1$.
Using Mehler's formula \cite{Bat53}
$$
\label{eqn:identity}
\exp \left( -\frac{(y-x z)^2}{1-z^2} \right)
= e^{-y^2} \sqrt{\pi(1-z^2)} 
\ \sum_{k=0}^{\infty}
\frac{z^k}{h_{k}} H_{k}(x) H_{k}(y),
$$
we will have the following expansions using
the Hermite polynomials.
For $1 \leq m < n \leq M+1$,
\begin{eqnarray}
\label{eqn:p_mn}
p_{t_n-t_m}(x, y) &=&
\frac{\sqrt{T}e^{-\frac{1}{2}(1+\gamma_{m}) (x/c_{m})^2}
e^{-\frac{1}{2}(1-\gamma_{n}) (y/c_{n})^2}}
{\sqrt{t_n(2T-t_m)}}
\\
&&\times \sum_{k=0}^{\infty}
\frac{e^{-k(\tau^{(n)}-\tau^{(m)})}}{h_{k}}
H_{k}\left( \frac{x}{c_{m}} \right) 
H_{k}\left( \frac{y}{c_{n}} \right),
\nonumber
\end{eqnarray}
and for $1 < m \leq M+1$,
\begin{eqnarray}
p_{t_1}(0,x)p_{t_m -t_1}(x,y)
&=& 
\frac{\sqrt{T}e^{- (x/c_1)^2}
e^{-\frac{1}{2}(1-\gamma_{m}) (y/c_{m})^2}}
{\sqrt{2\pi t_1 t_m (2T-t_1)}}
\nonumber\\
&&\times \sum_{k=0}^{\infty}
\frac{e^{-k(\tau^{(m)}-\tau^{(1)})}}{h_{k}}
H_{k}\left( \frac{x}{c_{1}} \right) 
H_{k}\left( \frac{y}{c_{m}} \right).
\nonumber
\label{eqn:p_1m}
\end{eqnarray}
Then from (\ref{def:Rk}), (\ref{def:Rkm})
and the orthogonal relation of the Hermitian polynomials,
we obtain
\begin{equation}
\label{eqn:Rkm}
R_{k}^{(m)}(x)
= \frac{e^{-\frac{1}{2} (1-\gamma_{m}) (x/c_{m})^2}
e^{k \tau^{(1)}}}{\sqrt{2\pi t_m}} 
\sum_{j=0}^{k} \alpha_{kj}
e^{-j \tau^{(m)}}
H_{j}\left(\frac{x}{c_{m}}\right).
\end{equation}
From the definition (\ref{def:Fmn}) and the expansion (\ref{eqn:p_mn})
we can obtain 
\begin{eqnarray}
\label{eqn:Fmn}
&& F^{m,n}(x,y) =
\frac{e^{-\frac{1}{2}(1+\gamma_{m})(x/c_m)^2} 
e^{-\frac{1}{2}(1+\gamma_{n})(y/c_n)^2} }{\sqrt{(2T-t_m)(2T-t_n)}} 
\\
&& \quad \times
\sum_{k=0}^{\infty} \sum_{\ell=0}^{\infty} 
\frac{ e^{k\tau^{(m)}}e^{\ell \tau^{(n)}}}{h_{k} h_{\ell}}
H_{k}\left(\frac{x}{c_{m}}\right) H_{\ell}\left(\frac{y}{c_{n}}\right)
\left\langle H_{k}\left(\frac{\cdot}{\sqrt{T}}\right), 
H_{\ell}\left(\frac{\cdot}{\sqrt{T}}\right) \right\rangle_{*},
\nonumber
\end{eqnarray}
where 
$\langle \cdot,\cdot \rangle_{*}$ is the antisymmetric 
inner product defined by
$$
\langle f, g \rangle_{*}=
\int_{-\infty}^{\infty} dw \int_{-\infty}^{w} dz \
e^{-(z^2 + w^2)/2T }
\Big[ f(z) g(w) - f(w) g(z) \Big].
$$
Put 
$R_{k}^{*}(x)=\sum_{j=0}^{k} \alpha_{kj} H_{j}
\left(\frac{x}{c_{M+1}}\right)$.
Then $\{R_{k}^{*}(x)\}$ satisfy the following
skew orthogonal relations
\begin{eqnarray}
&&\langle R^{*}_{2j}, R^{*}_{2\ell+1} \rangle_{*}
= - \langle R^{*}_{2\ell+1}, R^{*}_{2j} \rangle_{*}
=r^{*}_{j} \delta_{j \ell}, 
\label{eqn:skew03}\\
&&\langle R^{*}_{2j}, R^{*}_{2\ell} \rangle_{*}=0,
\quad
\langle R^{*}_{2j+1}, R^{*}_{2\ell+1} \rangle_{*}=0,
\quad \mbox{for} \ j, \ell=0,1,2, \dots,
\nonumber
\end{eqnarray}
where $r_{\ell}^{*} = 4h_{2\ell}T(c_1/2)^{4\ell+1}$.
We put 
$$
\beta_{k \, j} =\left\{
   \begin{array}{ll}
   2^{k} c_{1}^{-k} \delta_{j \, k},
   & \quad \mbox{if} \quad \mbox{ $k$ is even,}
   \\
  2^{k}\Big(\frac{k-1}{2} \Big)! 
  \Big\{ c_{1}^{j} (\frac{j-1}{2})! \Big\}^{-1},
   & \quad \mbox{if} \quad \mbox{$k,j$ are odd and $k \geq j$,} 
       \\
       0,
       & \quad \mbox{otherwise,} \\
   \end{array}\right.
$$
for nonnegative integers $k$ and $j$.
Then 
$\sum_{j=s}^{k} \beta_{kj} \alpha_{j s}=\delta_{ks}$,
if $0 \leq s \leq k$, and
\begin{equation}
H_{k}\left(\frac{x}{\sqrt{T}}\right)
=\sum_{j=0}^{k} \beta_{kj} R_{j}^{*}(x).
\label{eqn:H=betaR}
\end{equation}
From the definition (\ref{def:Phi})
and the equations (\ref{eqn:Rkm}) and (\ref{eqn:Fmn})
we have
\begin{eqnarray}
\Phi_{k}^{(m)}(x) 
&=& \frac{c_m}{\sqrt{2\pi t_m}\ (2T-t_m)}
e^{-\frac{1}{2}(1+\gamma_{m})(x/c_{m})^2}
e^{k \tau^{(1)}}  
\nonumber\\
&&\quad \times \sum_{\ell=0}^{\infty}\sum_{j=0}^{k} 
\frac{e^{ \ell \tau^{(m)}}}{h_{\ell}} 
H_{\ell}\left(\frac{x}{c_{m}} \right) \alpha_{kj} 
\left\langle H_{j} \left(\frac{\cdot}{\sqrt{T}} \right), 
H_{\ell}\left(\frac{\cdot}{\sqrt{T}} \right) \right\rangle_{*} 
\nonumber\\
&=&
\frac{e^{-\frac{1}{2}(1+\gamma_{m})(x/c_{m})^2}}
{\sqrt{2\pi T(2T-t_m)}}
e^{k \tau^{(1)}} 
\sum_{j=0}^{\infty} 
\langle R_{k}^{*}, R_{j}^{*} \rangle_{*}
\sum_{\ell=j}^{\infty}
\frac{e^{\ell \tau^{(m)}}}{h_{\ell}} H_{\ell}
\left( \frac{x}{c_{m}} \right) 
\beta_{\ell j}.
\nonumber
\end{eqnarray}
Using the skew orthogonal relations (\ref{eqn:skew03}),
we show that for $k=0,1,2, \dots$
\begin{eqnarray}
\label{eqn:Phi_even}
\Phi^{(m)}_{2k}(x) &=&
\frac{e^{-\frac{1}{2}(1+\gamma_{m})(x/c_{m})^2}r^{*}_{k}}
{\sqrt{2\pi T(2T- t_m)}}
  e^{2 k \tau^{(1)}}}
  \displaystyle{\sum_{\ell=2 k+1}^{\infty} 
  \frac{e^{\ell \tau^{(m)}}}{h_{\ell}}
  \beta_{\ell 2k+1} H_{\ell}
  \left( \frac{x}{c_{m}} \right), \\
\label{eqn:Phi_odd}
\Phi^{(m)}_{2k+1}(x) &=&
  - \frac{e^{-\frac{1}{2}(1+\gamma_{m})(x/c_{m})^2}r^{*}_{k}}
{\sqrt{2\pi T(2T- t_m)}}
  e^{(2 k+1) \tau^{(1)}}}
  \displaystyle{\sum_{\ell=2 k}^{\infty} 
  \frac{e^{\ell \tau^{(m)}}}{h_{\ell}}
  \beta_{\ell 2k} H_{\ell}
  \left( \frac{x}{c_{m}} \right).
\end{eqnarray}

Using above expansions we show the following lemma.
\vskip 3mm

\noindent{\bf Lemma 4.}\
{\itshape
For $1 \leq m, n \leq M+1$,
\begin{eqnarray}
&&F^{m,n}(x,y)=\sum_{k=0}^{\infty} \frac{1}{r_{k}}
\left[ \Phi_{2k}^{(m)}(x) \Phi_{2k+1}^{(n)}(y)
-\Phi_{2k+1}^{(m)}(x) \Phi_{2k}^{(n)}(y) \right],
\label{eqn:Fmn2}
\\
&&\widetilde{I}^{m,n}(x,y)=\sum_{k=N/2}^{\infty} 
\frac{1}{r_{k}}
\left[
\Phi_{2k}^{(m)}(x) \Phi_{2k+1}^{(n)}(y)
-\Phi_{2k+1}^{(m)}(x) \Phi_{2k}^{(n)}(y) \right].
\label{eqn:tImn}
\end{eqnarray}
}

\noindent{\it Proof.} \
By (\ref{eqn:Phi_even}), (\ref{eqn:Phi_odd}) and the relation
\begin{equation}
\label{eqn:r_k}
r_k = \frac{1}{2\pi T} \left(\frac{t_1}{2T-t_1} \right)^{2k+1/2} r_k^*,
\end{equation}
we have
\begin{eqnarray}
&&\sum_{k=0}^{\infty} \frac{1}{r_{k}}
\left[ - \Phi_{2k+1}^{(m)}(x) \Phi_{2k}^{(n)}(y)
+ \Phi_{2k}^{(m)}(x) \Phi_{2k+1}^{(n)}(y) \right]
\nonumber\\
&=& \frac{e^{-\frac{1}{2}(1+\gamma_{m})(x/c_{m})^2}
e^{-\frac{1}{2}(1+\gamma_{n})(y/c_{n})^2} }{\sqrt{(2T-t_m)(2T-t_n)}}
\sum_{k=0}^{\infty} r_{k}^{*}
\nonumber\\
&& \hskip 1.5cm \times \left\{ \sum_{j=2k}^{\infty} 
\frac{e^{j \tau^{(m)}}}{h_{j}} 
\beta_{j \, 2k} H_{j} \left( \frac{x}{c_{m}} \right) 
\times \sum_{\ell=2k+1}^{\infty} 
\frac{e^{\ell \tau^{(n)}}}{h_{\ell}}
\beta_{\ell \, 2k+1} 
H_{\ell} \left( \frac{y}{c_{n}} \right) \right. 
\nonumber\\
&& \hskip 2cm \left. - \sum_{j=2k+1}^{\infty} 
\frac{e^{j \tau^{(m)}}}{h_{j}} 
\beta_{j \, 2k+1} H_{j} \left(\frac{x}{c_{m}} \right) 
\times \sum_{\ell=2k}^{\infty} 
\frac{e^{\ell \tau^{(n)}}}{h_{\ell}}
\beta_{\ell \, 2k} H_{\ell} \left(\frac{y}{c_{n}} \right) \right\}. 
\nonumber
\end{eqnarray}
By (\ref{eqn:skew03}) and (\ref{eqn:H=betaR})
the right hand side of the above equation 
equals to 
\begin{eqnarray}
&& \frac{e^{-\frac{1}{2}(1+\gamma_{m})(x/c_{m})^2}
e^{-\frac{1}{2}(1+\gamma_{n})(y/c_{n})^2} }{\sqrt{(2T-t_m)(2T-t_n)}} 
\sum_{\mu=0}^{\infty} \sum_{\nu=0}^{\infty} 
\langle R_{\mu}^{*}, R_{\nu}^{*} \rangle_{*} 
\nonumber\\
&& \hskip 2cm \times
\sum_{j=\mu}^{\infty} \sum_{\ell=\nu}^{\infty} 
\frac{e^{j \tau^{(m)}} e^{\ell \tau^{(n)}}}{h_{j} h_{\ell}}
\beta_{j \mu} H_{j} \left( \frac{x}{c_{m}} \right) 
\beta_{\ell \nu} H_{\ell} \left( \frac{y}{c_{n}} \right) 
\nonumber\\
&=& F^{m,n}(x,y),
\nonumber
\end{eqnarray}
where we have used (\ref{eqn:Fmn}).
From the definitions (\ref{def:Imn}) and (\ref{def:tImn}), 
(\ref{eqn:tImn}) is derived from (\ref{eqn:Fmn2}).

\vskip 3mm

\SSC{Proof of Theorems}

The following formulae are known for (\ref{eqn:varphi})
\cite{Bat53,Sze75}.
For $u \in \R$,
\begin{eqnarray}
&&\lim_{\ell\to\infty}(-1)^{\ell} \ell^{1/4}\varphi_{2\ell}
\left(\frac{u}{2\sqrt{\ell}}\right)
= \frac{1}{\sqrt{\pi}}\cos u,
\label{eqn:cos}
\\
&&\lim_{\ell\to\infty}(-1)^{\ell} \ell^{1/4}\varphi_{2\ell+1}
\left(\frac{u}{2\sqrt{\ell}}\right)
= \frac{1}{\sqrt{\pi}}\sin u,
\label{eqn:sin}
\\
&&
\lim_{\ell\to\infty} 2^{-\frac{1}{4}}
\ell^{\frac{1}{12}} \varphi_{\ell}
\left( \sqrt{2\ell} - \frac{u}{\sqrt{2} \ \ell^{1/6}}\right) 
= {\rm Ai} (-u)
\label{eqn:Ai2}
\end{eqnarray}
Here we give the proof of Theorem 2 by using (\ref{eqn:Ai2}).
The proof of Theorem 1 will be easier and given
by the similar argument 
using (\ref{eqn:cos}) and (\ref{eqn:sin}).

Let 
$b^m(x)= \sqrt{2T-t_m} \exp \Big\{ 1/2 \gamma_m (x/c_m)^2 
- N \tau^{(m)} \Big\}$
and 
$\zeta^m(x)$ be the quaternion with 
$$
C(\zeta^m (x)) = \left( 
\begin{array}{cc}
b^m(x) & 0 \cr
\cr 0 & 1/b^m(x)
\end{array}
\right).
$$
For $\x^{(m)}_N \in \RV, 1 \leq m \leq M+1$,
we consider the transformation of the quaternions
$q^{m,n} \left(x_{i}^{(m)}, x_{j}^{(n)} \right) 
\mapsto \hq^{m,n} \left(x_{i}^{(m)}, x_{j}^{(n)} \right)$
defined by
$$
\hq^{m,n}(x,y)= \zeta^m(x) {q}^{m,n}(x,y) \zeta^n (y)^{-1}.
$$
We denote by 
$\hQ \left(\x^{(1)}_{N_1},\x^{(2)}_{N_2},\dots, 
\x^{(M+1)}_{N_{M+1}} \right)$
the self-dual
$\sum_{m=1}^{M+1} N_m \times \sum_{m=1}^{M+1} N_m$ 
quaternion matrix whose elements are
$\hq^{m,n} \left(x^{(m)}_i, x^{(n)}_j \right)$,
$1\le i\le N_m$, $1\le j\le N_n$,
$1\le m,n \le M+1$.
By the definition of quaternion determinants,
the following invariance is established:
$$
{\rm Tdet}\ Q \left(\x^{(1)}_{N_1},\x^{(2)}_{N_{2}},
\dots, \x^{(M+1)}_{N_{N_{M+1}}} \right)
= {\rm Tdet}\ \hQ \left(\x^{(1)}_{N_1},\x^{(m)}_{N_{2}},
\dots,\x^{(M+1)}_{N_{N_{M+1}}} \right).
$$
Hence to prove Theorem 2 it is enough to show the following lemma.

\vskip 3mm
\noindent{\bf Lemma 5.}\
{\itshape
Let $T_N= 2N^{1/3}$ and $t_m = T_N + s_m$, $1\leq m,n \leq M+1$.
Then for any $x, y \in \R$,
\begin{eqnarray}
&&\lim_{N\to\infty} \frac{1}{b^m(x)b^n(y)}
D^{m,n}(a_N(s_m) +x,  a_N(s_n)+y)
=\cD (s_m, x; s_n, y),
\label{eqn:D->cD}
\\
&&\lim_{N\to\infty} b^{m}(x)b^{n}(y)
\tI^{m,n}(a_N(s_m) +x, a_N(s_n)+y)
=\widetilde{\cI} (s_m, x; s_n, y),
\label{eqn:tI->tcI}
\\
&&\lim_{N\to\infty} \frac{b^{m}(x)}{b^{n}(y)}
\tS^{m,n}(a_N(s_m) +x, a_N(s_n)+y)
=\widetilde{\cS} (s_m, x; s_n, y).
\label{eqn:tS->tcS}
\end{eqnarray}
}

We start to prove this lemma by showing 
\begin{equation}
\label{eqn:p->cG}
\lim_{N\to\infty} \frac{b^{m}(x)}{b^{n}(y)}
p_{t_n-t_m}(a_N(s_m) +x, a_N(s_n)+y)
=\cG (s_m, x; s_n, y).
\end{equation}
By (\ref{eqn:p_mn}) and the fact
\begin{eqnarray}
&&\frac{a_N(s_m) +x}{c_m}= \sqrt{2N}+\frac{x}{\sqrt{2} N^{1/6}}
+ \cO(T_N^{-1}), 
\nonumber\\
&&\tau^{(n)} = \frac{s_n}{T_N} + \cO(T_N^{-2})
\nonumber
\end{eqnarray}
for large $N$, we have
\begin{eqnarray}
&&\lim_{N\to\infty} \frac{b^{m}(x)}{b^{n}(y)}
p_{t_n-t_m}(a_N(s_m) +x, a_N(s_n)+y)
\nonumber
\\
&&=\lim_{N\to\infty} 
\sqrt{\frac{1}{T_N}}
\sum_{p=-\infty}^{N} e^{\frac{p}{2N^{1/3}}(s_n-s_m)}
\varphi_{N-p} \left(\sqrt{2N}+\frac{x}{\sqrt{2} N^{1/6}}\right)
\varphi_{N-p}
\left(\sqrt{2N}+\frac{y}{\sqrt{2} N^{1/6}} \right)
\nonumber
\\
&&=\lim_{N\to\infty} 
\frac{1}{N^{1/3}}
\sum_{p=-\infty}^{N} e^{\frac{p}{2N^{1/3}}(s_n-s_m)}
{\rm Ai} \left(x +\frac{p}{N^{1/3}} \right)
{\rm Ai} \left(y +\frac{p}{N^{1/3}} \right),
\nonumber
\end{eqnarray}
where we have used (\ref{eqn:Ai2}).
Then we have (\ref{eqn:p->cG}).

From (\ref{def:Smn}), (\ref{eqn:Rkm}), (\ref{eqn:Phi_even})
and (\ref{eqn:Phi_odd}), we have
$$
S^{m,n}(x,y) = S_{1}^{m,n}(x,y)+S_{2}^{m,n}(x,y),
$$
with
\begin{eqnarray}
S_{1}^{m,n}(x,y) &=& \frac{b^{n}(y)}{c_n b^{m}(x)} 
\sum_{\ell=0}^{N-1} 
e^{(N- \ell)(\tau^{(n)}- \tau^{(m)})}
\varphi_{\ell}\left(\frac{x}{c_{m}} \right) 
\varphi_{\ell}\left(\frac{y}{c_{n}} \right), 
\nonumber\\
S_{2}^{(mn)}(x,y) &=& \frac{b^{n}(y)}{c_n b^{m}(x)} \varphi_{N-1}(y/c_{n})
\nonumber\\
&&\times
\sum_{k= N/2}^{\infty} \frac{B(N/2 + k)}{B(N/2 -1)}
e^{-(N- 2k-1)\tau^{(m)}}
\varphi_{2k+1}\left( \frac{x}{c_{m}} \right),
\nonumber
\end{eqnarray}
where
$B(k) = \frac{2^{k} k!}{\sqrt{ (2k+1)! } }$.
Since
$$
\frac{B(k)}{B(\ell)} = \left( \frac{k}{\ell} \right)^{1/4} 
\left(1+\cO\left( \frac{|k-\ell|}{k+\ell}\right)\right),
$$
by the same argument to show (\ref{eqn:p->cG}) we have 
\begin{equation}
\lim_{N\to\infty} \frac{b^{m}(x)}{b^{n}(y)}
S^{m,n}(a_N(s_m) +x, a_N(s_n)+y)
=\cS (s_m, x; s_n, y).
\label{eqn:S->cS}
\end{equation}
(\ref{eqn:tS->tcS}) is derived from (\ref{eqn:p->cG}) 
and (\ref{eqn:S->cS}).

From (\ref{eqn:Phi_even}) and (\ref{eqn:Phi_odd}),
by calculations with (\ref{eqn:r_k}), we have
\begin{eqnarray}
&&b^{m}(x)b^{n}(y)\Phi^{(m)}_{N+2p}(x)\Phi^{(n)}_{N+2p+1}(y)
\nonumber
\\
&=& -2^{3/2} r_{N/2+p} \frac{T}{\sqrt{N+2p+1}}
e^{(2p)\tau^{(n)}}
\varphi_{N+2p}\left( \frac{y}{c_n} \right)
\nonumber
\\
&&\qquad\qquad\times \sum_{k=p}^{\infty}
\frac{B(N/2 +k )}{B(N/2 +p)}
e^{(2k+1)\tau^{(m)}}
\varphi_{N+2k+1}\left( \frac{x}{c_m} \right).
\nonumber
\end{eqnarray}
From (\ref{eqn:tImn}) we obtain (\ref{eqn:tI->tcI})
by the same procedure as above.

From (\ref{eqn:Rkm}),
by calculations with (\ref{def:alpha}) and
$$
e^{-y^2/2}H_{\ell+1}(y) 
= -2 \frac{d}{dy}\left(e^{-y^2/2}H_{\ell}(y) \right)
+ 2\ell e^{-y^2/2}H_{\ell-1}(y),
$$
we have
\begin{eqnarray}
&&\frac{R^{(m)}_{2k}(x) R^{(n)}_{2k+1}(y)}{r_k b^{m}(x) b^{n}(y)}
\nonumber
\\
&=& \frac{-1}{2\sqrt{t_m t_n(2T-t_m)(2T-t_n)}}
e^{(N-2k)\tau^{(m)}}\varphi_{2k} \left(\frac{x}{c_m} \right)
e^{(N-2k+1)\tau^{(n)}}
\nonumber\\
&&\times \left\{ 
\varphi_{2k}^{\prime} \left(\frac{y}{c_n} \right)
+\sqrt{2k} \left(1- e^{-2\tau^{(n)}} \right)
\varphi_{2k-1}\left(\frac{y}{c_n} \right) \right\}.
\nonumber
\end{eqnarray}
Using the fact that
$$
\frac{(N-p)^{1/12}}{2^{3/4}N^{1/6}}
\varphi^{\prime}_{N-p} \left( \frac{a_N(s_n) +y}{c_n} \right)
= \left. \frac{d}{d\lambda}{\rm Ai}(y+\lambda)
\right|_{\lambda=p/N^{1/3}}+o(1)
$$
for large $N$, we obtain (\ref{eqn:D->cD}).
This completes the proof of Lemma 5.
\vskip 3mm

\footnotesize 

\vskip 1cm

\noindent 
Makoto Katori\\
Department of Physics\\
Chuo University \\
Kasuga, Bunkyo-ku, Tokyo 112-8551\\
Japan\\
e-mail: katori@phys.chuo-u.ac.jp

\vskip 1cm

\noindent 
Taro Nagao\\
Department of Physics\\
Graduate School of Science\\
Osaka University\\
Toyonaka, Osaka 560-0043\\
Japan\\
e-mail : nagao@sphinx.phys.sci.osaka-u.ac.jp

\vskip 1cm

\noindent 
Hideki Tanemura\\
Department of Mathematics and Informatics\\
Chiba University\\
1-33 Yayoi-cho, Inage-ku, Chiba 263-8522\\ 
Japan\\
e-mail : tanemura@math.s.chiba-u.ac.jp

\end{document}